\documentclass[11pt]{article}
\usepackage{amsfonts,latexsym,amsmath,amscd}
\usepackage[pdftex]{graphicx}
\newcommand{\cleqn}{\setcounter{equation}{0}}
\newcommand{\clth}{\setcounter{theorem}{0}}
\newcommand {\sectionnew}[1]{\section{#1}\cleqn\clth}

\newcommand{\beq}{\begin{equation}}
\newcommand{\eeq}{\end{equation}}
\newcommand{\beqa}{\begin{eqnarray}}
\newcommand{\eeqa}{\end{eqnarray}}
\newcommand{\beaa}{\begin{eqnarray*}}
\newcommand{\ben}{\begin{eqnarray*}}
\newcommand{\eaa}{\end{eqnarray*}}
\newcommand{\een}{\end{eqnarray*}}

\newcommand \nc {\newcommand}
\nc \proof {\noindent {\em{Proof.\/ }}} \nc \qed {$\Box$\hfill}
\newtheorem{theorem}{Theorem}[section]
\newtheorem{lemma}[theorem]{Lemma}
\newtheorem{proposition}[theorem]{Proposition}
\newtheorem{corollary}[theorem]{Corollary}
\newtheorem{definition}[theorem]{Definition}
\newtheorem{example}[theorem]{Example}
\newtheorem{remark}[theorem]{Remark}
\newtheorem{conjecture}[theorem]{Conjecture}
\newtheorem{question}[theorem]{Question}

\nc \bth[1] { \begin{theorem}\label{t#1} } \nc \ble[1] {
\begin{lemma}\label{l#1} } \nc \bpr[1] {
\begin{proposition}\label{p#1} } \nc \bco[1] {
\begin{corollary}\label{c#1} } \nc \bde[1] {
\begin{definition}\label{d#1}\rm } \nc \bex[1] {
\begin{example}\label{e#1}\rm } \nc \bre[1] {
\begin{remark}\label{r#1}\rm } \nc \bcon[1] {
\begin{conjecture}\label{con#1}\rm } \nc \bque[1] {
\begin{question}\label{que#1}\rm }
\nc {\eth} { \end{theorem} } \nc {\ele} { \end{lemma} } \nc {\epr} {
\end{proposition} } \nc {\eco} { \end{corollary} } \nc {\ede} {
\end{definition} } \nc {\eex} { \end{example} } \nc {\ere} {
\end{remark} } \nc {\econ} { \end{conjecture} } \nc {\eque} {
\end{question} }
\nc \thref[1]{Theorem \ref{t#1}} \nc \leref[1]{Lemma \ref{l#1}} \nc
\prref[1]{Proposition \ref{p#1}} \nc \coref[1]{Corollary \ref{c#1}}
\nc \deref[1]{Definition \ref{d#1}} \nc \exref[1]{Example \ref{e#1}}
\nc \reref[1]{Remark \ref{r#1}}

\def \Rset {{\mathbb R}}

\def \Gr { {\mathrm{Gr}} }

\nc \GRA[1] { \Gr_A^{(#1)} }   
\nc \GRAN { \GRA{N} } \nc \GrA[1] { \Gr_A(#1) }\nc \GrAa {
\GrA{\alpha} }
\nc \GRB[1] { \Gr_B^{(#1)} }   
\nc \GRBN { \GRB{N} } \nc \GrB[1] { \Gr_B(#1) } \nc \GrBb {
\GrB{\beta} }
\nc \GRMB[1] { \Gr_{MB}^{(#1)} }   
\nc \GRMBN { \GRMB{N} } \nc \GrMB[1] { \Gr_{MB}(#1) } \nc \GrMBb {
\GrMB{\beta} }

\begin{document}
\title{{\Large\bf{Maximal and minimal solutions of an Aronsson equation: $L^{\infty}$ variational problems versus the game theory.}}
\footnotetext{The author was partially supported by NSF grant
D0848378.}
\footnotetext{email: yyu1@math.uci.edu } }
\author{Yifeng Yu\\[2mm]
        Department of math, University of California at Irvine}

\date{}

\maketitle

\begin{abstract}

The Dirichlet problem
$$
\begin{cases}
\Delta_{\infty}u-|Du|^2=0  \quad \text{ on $\Omega\subset \Rset ^n$}\\
u|_{\partial \Omega}=g
\end{cases}
$$
might have many solutions, where $\Delta_{\infty}u=\sum_{1\leq
i,j\leq n}u_{x_i}u_{x_j}u_{x_ix_j}$. In this paper, we prove that
the maximal solution is the unique absolute minimizer for
$H(p,z)={1\over 2}|p|^2-z$ from calculus of variations in
$L^{\infty}$ and the minimal solution is the continuum value
function from the ``tug-of-war" game. We will also characterize
graphes of solutions which are neither an absolute minimizer nor a
value function. A remaining interesting question is how to interpret
those intermediate solutions. Most of our approaches are based on an
idea of Barles-Busca [BB].
\end{abstract}

\sectionnew{Introduction}

Let $\Omega$ be a bounded open set in $\Rset ^n$. Peres, Schramm,
Sheffield and Wilson in [PSSW] introduced a two-person differential
game called ``tug-of-war". Starting from a point $x\in \Omega$, at
each step with fixed length, two players toss a fair coin to
determine the order of move. One player tries to maximize the payoff
function and the other wants to minimize it. The game will stop if
one of them reaches the boundary of $\Omega$. In this paper, let us
 assume that the running payoff function is a constant
$-\tau$ and the terminal payoff function is $g\in
W^{1,\infty}(\Omega)$. Owing to [PSSW], as the step size tends to
zero, value functions of the game will converge to the unique
viscosity solution of the equation \beq{}
\begin{cases}{\Delta_{\infty}u\over
|Du|^2}=\tau \quad \text{on $\Omega$} \\

u|_{\partial \Omega}=g.
\end{cases}
\eeq Following the terminology in [PSSW], we call a viscosity
solution of equation (1.1) {\it a continuum value function of the
``tug-of-war" game}. See the {\it User's Guide} Crandall-Ishii-Lions
[CIL] for definitions of viscosity solutions of general nonlinear
elliptic equations. Here we should be careful about the operator
${\Delta_{\infty}u\over |Du|^2}$ when $|Du|$ vanishes. According to
 the definition in [PSSW], if a $C^2$ test function $\phi$ touches $u$ at
$x\in \Omega$ from above (or below) and $D\phi (x)=0$, we require
that $\max_{\{|p|=1\}}p\cdot D^2\phi (x)\cdot p\geq \tau$ (or
$\min_{\{|p|=1\}}p\cdot D^2\phi (x)\cdot p\leq \tau$). When $n=1$,

 $$
 \max_{\{|p|=1\}}p\cdot \phi'' (x)\cdot p=\min_{\{|p|=1\}}p\cdot \phi'' (x)\cdot
 p=\phi''(x).
 $$ Hence equation (1.1) is just $u''=\tau$.

Multiplying $|Du|^2$ on both side, we derive that the value function
is also a viscosity solution of the equation  \beq{}
\Delta_{\infty}u-\tau |Du|^2=0. \eeq However, except when $\tau=0$,
solutions of equation (1.2) might not be solutions of equation
(1.1). Here is a simple example.

\vspace{2mm}

{\bf Example I}: {\it $u_1=0$ and $u_2={1\over 2}x^2-{1\over 2}$ are
both smooth solutions of
$$
\begin{cases}
(u')^2u''-|u'|^2=0 \quad \text{on $(-1,1)$}\\
u(-1)=u(1)=0.
\end{cases}
$$
But $u_1=0$ is not a solution of $u''={(u')^2u''\over (u')^2}=1$}.

According to [PSSW], equation (1.1) admits a unique solution. But
the above example suggests that equation (1.2) might have multiple
solutions with prescribed boundary value. Equation (1.2) is a
so-called Aronsson equation associated to $H(p,z)={1\over
2}|p|^2-\tau z$. For general $H=H(p,z,x)\in C^1(\Rset ^n\times \Rset
\times \Omega)$, the correspondent Aronsson equation is
$$
A_{H}(u)=H_{p}(x,u,Du)\cdot D_{x}(H(x,u,Du))=0 \quad \text{in
$\Omega$}.
$$
Here $H_{p}$ is the partial derivative of $H$ with respect to $p$
and $D_{x}$ represents the derivative with respect to $x$ of
$H(x,u(x),Du(x))$. Aronsson equations are Euler-Lagrangian equations
for ``calculus of variations in $L^{\infty}$" which were initiated
by G. Aronsson in 60's ([A1-4]). Here is the general definition of
minimizers of such highly nonconventional variational problems. For
$H=H(p,z,x)\in C(\Rset ^n\times \Rset \times \Omega)$, we say that
$u\in W_{loc}^{1,\infty}(\Omega)$ is an {\it absolute minimizer for
$H$ in $\Omega$} if for any open set $V\subset \bar V\subset \Omega$
and $v\in W^{1,\infty}(V)$,
$$
u|_{\partial V}=v|_{\partial V}
$$
implies that
$$
\mathrm{esssup_{ V}}H(Du,u,x)\leq \mathrm{esssup_{V}}H(Dv,v,x).
$$
Crandall proved in [C] (see also Barron-Jensen-Wang [BJW]) that if
$H\in C^2$ and is quasiconvex in $p$, then an absolute minimizer for
$H=H(p,z,x)$ in $\Omega$ is a viscosity solution of the Aronsson
equation
$$
A_{H}(u)=0   \quad \text{in $\Omega$}.
$$
A function $f$ is {\it quasiconvex}  if the set $\{f<t\}$ are convex
for all $t\in \Rset$.

Let us focus on $H={1\over 2}|p|^2-\tau z$. Then any absolute
minimizer for H is a viscosity solution of equation (1.2) in
$\Omega$. However, except when $\tau=0$, the converse might not be
true. In
 Example I, $u_2={1\over 2}x^2-{1\over 2}$ is not an absolute minimizer.
In fact, $u_2|_{\partial \Omega}=0$, but
$$
{1\over 2}=\mathrm\mathrm{esssup_{(-1,1)}}({1\over
2}{(u_{2}^{'})}^2-u_2)>0.
$$

\vspace{2mm}

Hence two natural questions arise.

\vspace{2mm}

 (1) Is an absolute minimizer for $H$ unique with prescribed
boundary value?

\vspace{2mm}

(2) If uniqueness holds, what are the positions of the continuum
value function from the game theory and the absolute minimize among
all  viscosity solutions of equation (1.2)?

When $\tau=0$, equation (1.2) is the famous infinity Laplacian
equation. Jensen proved in [J] that Dirichlet problem of the
infinity Laplacian equation has a unique solution. Hence the
continuum value function and the absolute minimizer coincide in this
case. So let us look at $\tau\ne 0$. By properly scaling and
changing signs, we may assume that $\tau=1$. The following is our
main result. \bth{} Suppose that $g\in W^{1,\infty}(\Omega)$. Then
 there exists a unique absolute minimizer for
$H={1\over 2}|p|^2-z$ in $\Omega$ with boundary value $g$. The
absolute minimizer is the maximal viscosity solution of
 \beq{}\begin{cases}
\Delta_{\infty}u-|Du|^2=0  \quad \text{ in $\Omega$}\\
u=g \quad \text{ on $\partial \Omega$}.
\end{cases}
\eeq Moreover, the continuum value function  from the game theory is
the minimal viscosity solution of above equation.\eth

\bre{} In Example I, $u_1=0$ is the absolute minimizer and
$u_2={1\over 2}x^2-{1\over 2}$ is the continuum value function.
Also, it is easy to deduce from Theorem 1.1 that for general
$\tau>0$ ($\tau<0$), the absolute minimizer is the maximal (minimal)
solution and the continuum value function is the minimal (maximal)
solution. \ere

In Theorem 1.1, the uniqueness of an absolute minimizer follows
immediately after we prove that an absolute minimizer is the maximal
solution. There were various results on uniqueness of absolute
minimizers from $L^{\infty}$-variational problems. See for instance
Crandall-Gunnarsson-Wang [CGW], Jensen [J], Jensen-Wang-Yu [JWY],
Juutinen [Ju], Barles-Busca [BB], etc. However all those results
depend on uniqueness of solutions of Dirichlet problems for
correspondent Aronsson equations, which, as suggested by Example I,
might not hold in our case. To prove that an absolute minimizer is
the maximal solution, we first use an idea from [BB] to reduce
inhomogeneous boundary conditions to homogeneous boundary
conditions. Then, by combine use of the PDE (1.3) and the definition
of absolute minimizers,  we prove that if an absolute minimizer
vanishes on the boundary, then it must be zero.

 We want to point out that the existence of absolute minimizers does not follow directly
 from the usual $L^{p}$ approximation introduced by Aronsson (see [BJW]) since
 $H={1\over 2}|p|^2-z$ is not bounded from below. What we do is to introduce an
 auxiliary $\hat H\geq 0$ and show that
 absolute minimizers for $\hat H$ are also absolute minimizers for
 $H={1\over 2}|p|^2-z$. Our approach relies on the fact that any
 solution of equation (1.3) is bounded from above by its maximum value on
 $\partial \Omega$. This is because a viscosity solution of equation (1.3) is a viscosity subsolution of the infinity
 Laplacian equation.

 \vspace{2mm}

{\bf Outline of our paper}. In section 2, we will prove Theorem 1.1.
In section 3, we give a characterization of solutions of equation
(1.3) which are neither the absolute minimizer nor the continuum
value function. A remaining interesting question is how to interpret
those solutions. \vspace{2mm}

 {\bf Notations.} We denote $B_r(x_0)$ as an open ball centered at
 $x_0$ with radius $r$. For $\delta>0$, we write
 $$
 \Omega_{\delta}=\{x\in \Omega|\ d(x,\partial \Omega)>\delta\}.
 $$
If $V$ is a subset of $\Rset ^n$, $\partial V$ denotes its boundary
and $\bar V$ the closure.  Moreover, if $f$ is a semiconvex
function, i.e. $f(x)+C|x|^2$ is convex for some $C>0$, we denote
$D^{-}f(x_0)$ as the subdifferentials of $f$ at $x_0$. That is
 $$
 D^{-}f(x_0)=\{p\in \Rset^n|\ f(x)\geq f(x_0)+p\cdot
 (x-x_0)-o(|x-x_0|)\},
 $$
where $o(|x-x_0|)$ means that $\lim_{x\to x_0}{o(|x-x_0|)\over
|x-x_0|}=0$.

\bre{} From now on, $H={1\over 2}|p|^2-z$. Moreover, we use
``absolute minimizer(s)" as an abbreviation for ``absolute
minimizer(s) for $H$ in $\Omega$" unless we specify the functional.
\ere

\sectionnew{Proofs}

We first use an idea from [BB] to prove a key lemma.

\ble{} Suppose that $u\in C(\bar \Omega)$ is a semiconvex viscosity
subsolution of equation \beq{} \Delta_{\infty} u-|Du|^2=0  \quad
\text{in $\Omega$}\eeq and $v\in C(\bar \Omega)$ is a viscosity
solution of the above equation. Assume that $$\max_{\bar
\Omega}(u-v)>\max_{\partial \Omega}(u-v). $$If
$u(x_0)-v(x_0)=\max_{\bar \Omega}(u-v)$ for some $x_0\in \Omega$,
then there exists $r_0>0$ such that
$$
u(x)=u(x_0)  \quad \text{for $x\in B_{r_0}(x_0)$}.
$$
\ele \proof For $\delta>0$ and $h\in B_{\delta}(0)$, denote
$M_{\delta}(h)=\max_{\bar \Omega_\delta}(u(x+h)-v(x))$.  It is clear
that $M_{\delta}(h)$ is a semiconvex function of $h$. Since the
maximum value of $u-v$ is not attained $\partial \Omega$, there
should exist $\delta_1>0$ such that for all $h\in B_{\delta_1}(0)$,
\beq{} \{x\in \bar \Omega_{\delta_1}|\
u(x+h)-v(x)=M_{\delta_1}(h)\}\subset \Omega_{2\delta_1}. \eeq Now I
claim that \beq{} 0\in D^{-}M_{\delta_1}(h)  \quad \text{for all
$h\in B_{\delta_1}(0)$}. \eeq In fact, fix $h$ and let us denote
$$
w_{\epsilon,h}(x,y)=(1+\epsilon) u(x+h)-v(y)-{1\over
2\epsilon}|x-y|^2.
$$
Suppose that $(\bar x,\bar y)\in \{(x,y)\in \bar
\Omega_{\delta_1}\times \bar \Omega_{\delta_1}|\ w_{ \epsilon,
h}(\bar x,\bar y)=\max_{x,y\in \bar \Omega_{\delta_1}}w_{\epsilon,
h}\}$. Owing to (2.2), when $\epsilon$ is small enough, we have that
$(\bar x,\bar y)\in \Omega_{\delta_1}\times \Omega_{\delta_1}$.
According to
[CIL], there exist $X$ and $Y$ such that\\[2mm]

(1) $({\bar x-\bar y\over \epsilon},X)\in  {\bar
J}_{\Omega_{\delta_1}}^{2,+}\mathrm{[}(1+\epsilon)u(\bar
x+h)\mathrm{]},\quad ({\bar x-\bar y\over \epsilon},Y)\in {\bar
J}_{\Omega_{\delta_1}}^{2,-}v(\bar
y),$\\[2mm]

(2) $-{3\over \epsilon}I_n\leq X\leq Y\leq {3\over \epsilon}I_n$.

\vspace{2mm}

Here $X$, $Y$, $\bar x$ and $\bar y$ all depend on $\epsilon$. See
[CIL] for definitions of ${\bar J}_{V}^{2,+}$ and ${\bar
J}_{V}^{2,-}$. Owing to equation (2.1), we have that
$$
{\bar x-\bar y\over \epsilon}\cdot X\cdot {\bar x-\bar y\over
\epsilon}\geq (1+\epsilon)|{\bar x-\bar y\over \epsilon}|^2
$$
and
$$
{\bar x-\bar y\over \epsilon}\cdot Y\cdot {\bar x-\bar y\over
\epsilon}\leq |{\bar x-\bar y\over \epsilon}|^2.
$$
Due to (2) above, we must have that
$$
{\bar x-\bar y\over \epsilon}=0.
$$
Since $u$ is semiconvex, $u$ is differentiable at $\bar x+h$ and
$$
D^{-}u(\bar x+h)=\{0\}.
$$
Passing  to a subsequence if necessary, we may assume that
$$
\lim_{\epsilon\to 0}\bar x=\lim_{\epsilon\to 0}\bar y=z_0.
$$
It is clear that $z_0\in \{x\in \bar \Omega_{\delta_1}|\
u(x+h)-v(x)=M_{\delta_1}(h)\}$. Since $u(\cdot+h)$ is semiconvex,
the set $D^{-}u(x)$ is upper-semicontinuous. Therefore
$$
0\in D^{-}u(z_0+h).
$$
Hence
$$
u(z_0+\hat h)\geq u(z_0+h)-o(|\hat h-h|).
$$
Therefore
$$
M_{\delta_1}(\hat h)\geq u(z_0+\hat h)-v(z_0)\geq
u(z_0+h)-v(z_0)-o(|\hat h-h|)=M_{\delta_1}(h)-o(|\hat h-h|).
$$
So
$$
0\in D^{-}M_{\delta_1}(h).
$$
Hence our claim holds. Therefore
$$
M_{\delta_1}(h)=M_{\delta_1}(0)  \quad \text{for $|h|\leq
\delta_1$}.
$$
Accordingly,
$$
u(x_0+h)-v(x_0)\leq M_{\delta_1}(h)=M_{\delta_1}(0)=u(x_0)-v(x_0).
$$
This implies that
$$
u(x_0+h)\leq u(x_0)  \quad \text{for $|h|\leq \delta_1$}.
$$
Since $u$ is a viscosity subsolution of equation (2.1), $u$ is a
viscosity subsolution of the infinity Laplacian equation
$$
\Delta_{\infty}u=0.
$$
According to the well known differential Harnack inequality (see
Lemma 2.5 in [CEG] for instance), we must have that \beq{}
u(x_0+h)=u(x_0) \quad \text{for $|h|\leq \delta_1$}. \eeq \qed

The following lemma says that the graph of an absolute minimizer can
not contain wells. Its proof is a combine use of the PDE and the
definition of absolute minimizers. \ble{} Suppose that $V$ is a
bounded open set in $\Rset ^n$. Assume that $w$ is an absolute
minimizer for $H$ on $V$ and
$$
w=c  \quad \text{on $\partial V$}.
$$
Then
$$
w\equiv c  \quad \text{in $V$}.
$$
\ele \proof Since $w-c$ is also an absolute minimizer, we may assume
that $c=0$. Since $w$ is an absolute minimizer, it is a viscosity
solution of equation (2.1). So it is a viscosity subsolution of the
infinity Laplacian equation
$$
\Delta_{\infty}w=0.
$$
Owing to the {\it maximum principle} for the infinity Laplacian
equation, we have that
$$
w\leq 0  \quad \text{in $V$}.
$$
Since $w$ is an absolute minimizer and vanishes on the boundary,
according to the definition of absolute minimizers,
$$
\mathrm{esssup}_{x\in V}(|Dw|^2-w)\leq 0.
$$
So
$$
w\geq 0  \quad \text{in $V$}.
$$
Therefore,
$$
w\equiv 0.
$$
\qed

Next lemma says that graphs of continuum value functions can not
contain flat pieces.

\ble{} Suppose that $u$ is a viscosity subsolution of equation (1.1)
with $\tau=1$. Then there does not exist a nonempty open set $V$
such that
$$
u\equiv \mathrm{constant}  \quad \text{in $V$}.
$$
\ele \proof We argue by contradiction. Suppose that there exists
such $V$. Choose a point $x_0\in V$. Then the quadratic polynomial
$$
P(x)=u(x_0)+{1\over 4}|x-x_0|^2
$$
touches $u$ at $x_0$ from the above in $V$. Since $DP(x_0)=0$, owing
to the definition of viscosity subsolutions  of equation (2.1), we
should have that
$$
\max_{\{|\xi|=1\}}\xi\cdot D^2P(x_0)\cdot \xi\geq 1.
$$
However, $\max_{\{|\xi|=1\}}\xi\cdot D^2P(x_0)\cdot \xi={1\over 2}$.
This is a contradiction. Therefore our lemma holds.
\qed

\vspace{2mm}

{\bf Proof of Theorem 1.1.} {Step I: (Existence of an absolute
minimizer)}. We may assume that $g\leq 0$. Now let us consider a new
Hamiltonian
$$
\hat H(p,z)={1\over 2}|p|^2-z^{-},
$$
where $z^{-}=\min\{z,0\}$. Clearly, $\hat H\geq {1\over 2}|p|^2$.
So, the existence of an absolute minimizer for $\hat H$ with
boundary value $g$ follows from the usual $L^{p}$ approximation. See
for instance [BJW]. Suppose that $w$ is an absolute minimizer for
$\hat H$ with boundary value $g\leq 0$. I want to show that $w$ is
also an absolute minimizer for $H$. In fact, assume that $V$ is an
open subset of $\Omega$ and $f\in W^{1,\infty}(V)$ such that
$$
f=w  \quad \text{on $\partial V$}.
$$
We need to prove that  \beq{} \mathrm{esssup}_{V}({1\over
2}|Dw|^2-w)\leq \mathrm{esssup}_{V}({1\over 2}|Df|^2-f).\eeq First I
claim that $w\leq 0$. We argue by contradiction. If not, since
$w|_{\partial \Omega}\leq 0$, then there exists an open subset
$U\subset \bar U\subset \Omega$ such that
$$
w>0 \quad \text{in $U$}
$$
and
$$
w=0 \quad \text{on $\partial U$}.
$$
Since $w$ is an absolute minimizer for $\hat H$,
$$
\mathrm{esssup}_{U}\hat H(Dw,w)\leq \mathrm{esssup}_{U}\hat
H(0,0)=0.
$$
Since $\hat H\geq {1\over 2}|p|^2$, we get that
$$
|Dw|=0 \quad \text{a.e in $U$}.
$$
Accordingly, $w\equiv 0$ in $U$. This is a contradiction.  Therefore
our claim holds, i.e, $w\leq 0$ in $\Omega$. Hence
$$
f\leq 0  \quad \text{on $\partial V$}.
$$
Therefore,
$$
\mathrm{esssup}_{V}({1\over 2}|Df|^2-f)\geq 0.
$$
Hence  \beq{}
 \mathrm{esssup}_{V}({1\over
2}|Df|^2-f)\geq \mathrm{esssup}_{V}({1\over 2}|Df^{-}|^2-f^{-}).
\eeq Note that
$$
{1\over 2}|Df^{-}|^2-f^{-}=\hat H(Df^{-},f^{-}).
$$
Since $w$ is an absolute minimizer for $\hat H$ and $w=f=f^{-}$ on
$\partial V$, we have that \beq{} \mathrm{esssup}_{V}\hat
H(Dw,w)\leq \mathrm{esssup}_{V}\hat
H(Df^{-},f^{-})=\mathrm{esssup}_{V}({1\over 2}|Df^{-}|^2-f^{-}).\eeq
Since $w\leq 0$, $w=w^{-}$. Hence \beq{} {1\over 2}|Dw|^2-w=\hat
H(Dw,w) \quad \text{a.e in $\Omega$}. \eeq Combining (2.6)-(2.8) ,
(2.5) holds. So $w$ is indeed an absolute minimizer for $H={1\over
2}|p|^2-z$.

\vspace{2mm}

{Step II: Next we show that an absolute minimizer is the maximal
viscosity solution of equation (1.3). Assume that $w$ is an absolute
minimizer and $u$ is an arbitrary viscosity subsolution of equation
(1.3). Our goal is to prove that \beq{} w\geq u \quad \text{in $\bar
\Omega$}. \eeq By considering super-convolution of $u$ and routine
modifications, we may assume that $u$ is semiconvex. If (2.9) does
not hold, there must exist $x_0\in \Omega$ such that
$$
u(x_0)-w(x_0)=\max_{\bar \Omega}(u-w)>0.
$$
According to Lemma 2.1, there exists $r>0$ such that $\overline
B_r(x_0)\subset \Omega$ and
$$
u\equiv u(x_0)   \quad \text{in $B_{r}(x_0)$}.
$$
We say that $\mathcal {O}$ is an {\it admissible open subset of
$\Omega$} if $x_0\in \mathcal {O}$ and
$$
u\equiv u(x_0)\quad \text{in $\mathcal {O}$}.
$$
Denote
$$
V=\cup_{\{\text{$\mathcal {O}$ is an admissible open subset of
$\Omega$}\}}\mathcal {O}.
$$
Note that $V$ is not empty since $B_r(x_0)\subset V$. I claim that
for $y\in \partial V$,
$$
w(y)>w(x_0).
$$
Owing to the choice of $x_0$, it is clear that $w(y)\geq w(x_0)$. If
$y\in \partial \Omega$, it is easy to see that
$w(y)=u(y)=u(x_0)>w(x_0)$. If $y\in \Omega$ and $w(y)=w(x_0)$, then
$$
u(y)-w(y)=u(x_0)-w(x_0)=\max_{\bar \Omega}(u-w)>0.
$$
By Lemma 2.1, there exists $r'>0$ such that
$$
u\equiv u(y)=u(x_0)  \quad \text{$B_{r'}(y)$}.
$$
Hence $B_{r'}(y)\cup V$ is an admissible open subset of $\Omega$. By
the definition of $V$, we have that  $B_{r'}(y)\subset V$. This
contradicts to $y\in \partial V$. Hence our claim holds.
Accordingly, there must exist a $\delta>0$ and an open subset
$V'\subset \bar V'\subset V$ such that $x_0\in V'$, \beq{}
w(x)<w(x_0)+\delta  \quad \text{in $V'$}\eeq and

$$
w(x)=w(x_0)+\delta  \quad \text{on $\partial V'$}.
$$
Hence by Lemma 2.2, $w\equiv w(x_0)+\delta$ in $V'$. This
contradicts to (2.10). Therefore (2.9) holds.

\vspace{2mm}

{Step III:} Finally, we need to show that the continuum value
function from the ``tug-of-war" game is the minimal viscosity
solution of equation. Suppose $u$ is a viscosity subsolution of
equation (1.1) and $v$ is an arbitrary viscosity solution of
equation (1.3). We need to show that \beq{} v\geq u  \quad \text{in
$\Omega$}.\eeq By super-convolution and routine modifications, we
may assume that $u$ is semiconvex. We argue by contradiction. If
$(2.11)$ is not true, owing to Lemma 2.1, there must exist a
nonempty open subset $V$ of $\Omega$ such that
$$
u\equiv c  \quad \text{in $V$}.
$$
This is impossible according to Lemma 2.3. Hence (2.11) holds.
\qed

\vspace{2mm}

The following theorem provides an alternative way to see why the
 continuum value function is the minimal solution.

\bth{} Any viscosity solution $u$ of equation (1.3) is a viscosity
supersolution of equation (1.1) with $\tau=1$.\eth

 \proof We argue by contradiction.
If not, then there exists $x_0\in \Omega$ and $\phi \in C^2(\Omega)$
such that \beq{} \phi(x)-u(x)<\phi (x_0)-u(x_0)=0  \quad \text{for
$x\in \Omega\backslash \{x_0\}$} \eeq and
$$
\min_{\{|p|=1\}}p\cdot D^2\phi (x_0)\cdot p>1.
$$
Hence the least eigenvalue of the Hessian matrix $D^2\phi (x_0)$
must be larger than $1$. Therefore
 \beq{} D^2\phi (x_0)>I_n, \eeq
where $I_n$ is the $n\times n$ identity matrix. Since $u$ is a
viscosity solution of (1.3), we have that
$$
D\phi (x_0)\cdot D^2\phi (x_0)\cdot D\phi (x_0)\leq |D\phi (x_0)|^2.
$$
By (2.13), $D\phi (x_0)=0$. Also owing to (2.13), there exists
$\delta>0$ such that \beq{} D\phi (x)\ne 0  \quad \text{for $x\in
B_{\delta}(x_0)\backslash \{x_0\}$}. \eeq For $h\in B_{r}(0)$, we
choose $x_h\in \Omega_{r}$ such that
$$
\phi (x_h+h)-u(x_h)=\max_{\bar \Omega_{r}}(\phi (x+h)-u(x)).
$$
According to (2.12), it is easy to see that when $r$ is small, $x_h$
will be close to $x_0$. Hence when $r$ is small enough, we have that
\beq{} D^2\phi (x_h+h)>I_n. \eeq Since $u$ is a viscosity solution
of (1.3), we have that
$$
D\phi (x_h+h)\cdot D^2\phi (x_h+h)\cdot D\phi (x_h+h)\leq |D\phi
(x_h+h)|^2.
$$
According to (2.15),
$$
D\phi (x_h+h)=0.
$$
By (2.14), when $r$ is sufficiently small, we must have that
$$
x_h+h=x_0.
$$
Hence due to the choice of $x_h$,
$$
\phi (x_0)-u(x_0-h)\geq \phi (x_0+h)-u(x_0).
$$
So
$$
\phi (x_0+h)+\phi (x_0-h)\leq \phi (x_0+h)+u(x_0-h)\leq 2\phi (x_0).
$$
This contradicts to (2.13) when $h$ is small.  Hence our claim
holds. \qed

\sectionnew{Other solutions of equation (1.3)}

In this section, we will give a characterization of graphes of
 {\it intermediate solutions}, i.e. those solutions between the absolute minimizer and the continuum value
function. Before  stating the theorem, we define some terminologies.

We say that the graph of a function $f\in C(\bar \Omega)$ has {\it a
well} if there exists a open set $V\subset \Omega$ such that
$$
\min_{\bar V}f< \min_{\partial V}f.
$$

We say that the graph of $f\in C(\bar \Omega)$ has a flat piece if
$f$ is constant in some open subset of $\Omega$.

\bth{} Suppose that $u$ is a viscosity solution of equation (1.3).
Then

(i) $u$ is not the absolute minimizer if and only its graph has
wells.

(ii) $u$ is not the value function if and only if its graph has flat
pieces.

Especially, $u$ is an intermediate solution if and only if its graph
 has both wells and flat pieces.
 \eth

\proof  (i) Note that in Step II of the proof of Theorem 1.1, we
only use the fact that the graph of an absolute minimizer has no
well. Hence (i) holds.

(ii) The sufficiency part of (ii) is Lemma 2.3. Hence we only need
to prove the necessity part. Assume that $v$ is the viscosity
solution of equation (1.1) with $\tau=1$. Suppose that $u\ne v$. We
are going to show that there exists a open set $U\subset \Omega$
such that $u$ is constant in $U$. Since $u\ne v$, we have that
$$
\max_{\bar \Omega}(u-v)>0.
$$
Hence there must exist $\delta>0$ such that for $h\in B_{\delta}(0)$
$$
\{x\in \bar \Omega_{\delta}|\ u(x+h)-v(x)=\max_{\bar
\Omega_{\delta}}(u(\cdot+h)-v)\}\subset \Omega_{3\delta}.
$$
Now fix $\delta$. For $\epsilon>0$, we denote $u_{\epsilon}$ as the
super-convolution of $u$, i.e,
$$
u_{\epsilon}(x)=\max_{y\in \bar \Omega}( u(y)-{1\over
\epsilon}|x-y|^2).
$$
It is clear that when $\epsilon$ is small enough, $u_{\epsilon}$ is
a viscosity subsolution of equation (1.3) in $\Omega_{\delta\over
2}$ and for $h\in B_{\delta\over 4}(0)$
$$
\{x\in \bar \Omega_{\delta}|\ u_{\epsilon}(x+h)-v(x)=\max_{\bar
\Omega_{\delta}}(u_{\epsilon}(\cdot+h)-v)\}\subset \Omega_{2\delta}.
$$
Note that $u_{\epsilon}$ is semiconvex. Choose $x_{\epsilon}\in
\Omega_{\delta}$ such that
$$
u_{\epsilon}(x_{\epsilon})-v(x_{\epsilon})=\max_{\Omega_{\delta}}(u_{\epsilon}-v).
$$
Owing to (2.4),
$$
u_{\epsilon}\equiv u_{\epsilon}(x_{\epsilon})  \quad \text{ in
$B_{\delta\over 4}(x_{\epsilon})$}.
$$
Passing to a subsequence if necessary, we may assume that
$$
\lim_{\epsilon\to 0}x_{\epsilon}=x_0\in \bar \Omega_{\delta}.
$$
Then
$$
u\equiv u(x_0) \quad \text{ in $B_{\delta\over 4}(x_{0})$}.
$$
\qed

\bre{} Equation in Example I actually possesses infinitely many
intermediate solutions. This motivates us to ask the following two
questions which we will investigate in the future.\\[2mm]
Q1  Is it true that there are infinitely many intermediate solutions
if the absolute minimizer and the value function do not
coincide?\\[2mm]
Q2  How to interpret those intermediate solutions?\ere

\renewcommand{\em}{\textrm}
\begin{small}
  \renewcommand{\refname}{ {\flushleft\normalsize\bf{References}} }

\end{small}

\end{document}